\documentclass[12pt,a4paper]{amsart}

\usepackage{amsmath,amssymb,amscd,array, mathrsfs,hyperref,bbm,diagbox}

\usepackage{tabularx,mdwtab}

\usepackage{DLdef1}

\newcommand{\del}{\partial}
\rmnameii{etr}{etr}
\rmnameii{evv}{ev}

\begin{document}

\title{Computer-aided study of double extensions of restricted Lie superalgebras preserving the non-degenerate closed 2-forms in characteristic 2}

\author{Sofiane Bouarroudj$^1$, Dimitry Leites$^2$, Jin Shang$^1$}

\address{$^1$New York University Abu Dhabi, Division of Science and Mathematics, P.O. Box 129188, Abu Dhabi\\ United Arab Emirates;\\
 \{sb3922, js8544\}@nyu.edu;\\
 $^2$Department of Mathematics Stockholm University, Albanov\"agen 28, SE-114 19, Stockholm, Sweden;\\ dimleites@gmail.com}

% \thanks will become a 1st page footnote.
%\thanks{SB was supported by the grant NYUAD-065.}

% \thanks will become a 1st page footnote.
%\thanks{SB was supported by the grant NYUAD-065.}

% General info

%\date{\today}


\begin{abstract} A Lie (super)algebra with a non-degenerate invariant symmetric bilinear form $B$ is called a NIS-(super)algebra. The double extension $\mathfrak{g}$ of a NIS-(super)algebra $\mathfrak{a}$ is the result of simultaneous adding to $\mathfrak{a}$ a central element and a derivation so that $\mathfrak{g}$ is a NIS-algebra. Loop algebras with values in simple complex Lie algebras are most known among the Lie (super)algebras suitable to be doubly extended. In characteristic 2 the notion of double extension acquires specific features. 

Restricted Lie (super)algebras are among the most interesting modular Lie superalgebras. In characteristic 2, using Grozman's Mathematica-based package SuperLie, we list double extensions of restricted Lie superalgebras preserving the non-degenerate closed 2-forms with constant coefficients. The results are proved for the number of indeterminates ranging from 4 to 7 --- sufficient to conjecture the pattern for larger numbers. Considering multigradings allowed us to accelerate computations up to 100 times.
\end{abstract}


\makeatletter
\@namedef{subjclassname@2020}{\textup{2020} Mathematics Subject Classification}
\makeatother

\maketitle

\markboth{\itshape Sofiane Bouarroudj\textup{,} Dimitry Leites\textup{,} Jin Shang}{{\itshape Double extensions of Lie superalgebras}}

 \thispagestyle{empty}

\setcounter{tocdepth}{2}
\tableofcontents

\section {Introduction}

\subsection{Setting of the problem} For a given Lie (super)algebra $\mathfrak{a}$ over the ground field $\Kee$, the notion of double extension $\fg:=\mathscr{K} \oplus \fa \oplus \mathscr{K} ^*$, where $\mathscr{K} :=\Kee c$, and $\mathscr{K}^* :=\Kee D$, was recently identified, see \cite{MR}. It simultaneously involves 

1) a central extension of $\mathfrak{a}$ (with center $\mathscr{K} :=\Kee c$), 

2) a derivation $D$ of $\mathfrak{a}$, 

3) a non-degenerate invariant symmetric bilinear form (briefly: NIS) on $\mathfrak{a}$ extendable to $\fg$. Hereafter any Lie (super)algebra with a NIS is called a \textit{NIS-(super)algebra}.

Most known examples of double extensions are affine Kac--Moody algebras $\fg$ over $\Cee$ or $\Ree$ important in physics and mathematics (Google returns ca 500K entries); here $\mathfrak{a}:=\fs^{\ell(1)}$ is a loop algebra with values in a simple finite-dimensional Lie algebra $\fs$, see \cite{K}. Less known examples: the Lie superalgebra $\fgl(a|a)$ over $\Cee$ which is a double extension of $\mathfrak{psl}(a|a)$, and the Lie algebra $\fgl(np)$ in characteristic $p>0$ which is a double extension of $\mathfrak{psl}(np)$. 

In what follows we show that for the double extension to be ``interesting", i.e., not just the direct sum of two ideals ($\fa$ and $\mathscr{K} \oplus \mathscr{K} ^*$), the following conditions should be satisfied: 

a) the derivation $D$ of $\mathfrak{a}$ \textit{must be outer} for any $p$; if $p=2$, and $D\in(\mathfrak{out}\ \fa)_\od$, then condition $D^2=0$ is a must, see \cite{BKLS};

b) the central extension has to be non-trivial, see Subsection~\ref{nontrCE}.

The simple Lie (super)algebras with NIS are classified for various types of Lie (super)algeb\-ras over the field $\Kee$ of characteristic $p\neq 2$ in review \cite{BKLS}. 

The new and most interesting results of \cite{BeB1} are general constructions and examples of double extensions of Lie superalgebras for $p=2$. Here we recall the main definitions and general results of \cite{BeB1} carefully pointing at the difference between the cases where $p\neq 2$ and the cases where $p=2$. 

In particular, in \cite{BeB1}, the double extensions of simple NIS-Lie superalgebras $\fh_\Pi^{(1)}(0|4)$ and 
$\fh_\Pi^{(1)}(0|5)$, where $\fg^{(1)}$ is the derived algebra of $\fg$, see eq.~\eqref{prim}, were classified. Clearly, $\fh_\Pi^{(1)}(0|n)$ has the Poisson Lie superalgebra $\fpo_\Pi(0|n)$ as its double extension, but there are more: several new double extensions were found; one of them gave an interpretation of a result in \cite{BGLL1}. 

P.~Deligne advised one of the co-authors of this note to consider restricted algebras first of all, as pertaining to geometry, see \cite{LL}. In \cite{BBH}, the notion of double extension was extended to \textit{restricted} Lie algebras. 
If $p=2$, there are several notions of restrictedness, see \cite{BLLS}; here we consider the ``classical'' one.
Vectorial Lie (super)algebras can only be restricted if the shearing vector is equal to $\mathbbm{1}:=(1,\dots,1)$. 

\textbf{Hereafter $\Kee$ is an algebraically closed field of characteristic $2$ and the shearing vector is equal to $\mathbbm{1}$}, so we do not indicate it. For details of the description of simple Lie superalgebras we study in what follows, see \cite{BGLLS}.

\textbf{Our results}. There are two types of super analogs of the Hamiltonian Lie algebra: series $\fh$, and series $\fle$ introduced in \cite{Le}. If $\Char \Kee>0$, there are several non-isomorphic versions of these series; here we consider the ``standard'' ones, i.e., with constant coefficients, see \cite[Subsection 4.7.12]{BKLS}:

a) the restricted Lie (super)algebras $\mathfrak{h}_B^{(1)}(a |b)$  with forms $B=\Pi\Pi$, $\Pi I$, $I\Pi$ and $II$ on the superspace spanned by the $X_i$, where 
${X=(X_1,\dots,X_{a+b})}$ form $a$ even and $b$ odd indeterminates (in particular, $B=\Pi$ and $I$ on the space, no ``super''); 

b) the restricted Lie superalgebras $\mathfrak{le}^{(1)}(a |a)$ on $a$ even and $a$ odd indeterminates.

We sharpen the conjecture from \cite{BeB1} on the shape of double extensions $\mathfrak{h}_\Pi^{(1)}(0|n)$ for $n>5$ and prove it. We computed double extensions of $\mathfrak{h}_B^{(1)}(a |b)$ for $4\leq a+b\leq 6$ and $\mathfrak{h}_\Pi^{(1)}(0|7)$, and $\mathfrak{le}^{(1)}(a |a)$ for $a=2,3$: this suffices to see the pattern for any $a$ and $b$. 

The deforms of series $\fh$ and $\fle$ preserving the 2-forms with non-constant coefficients (not classified yet, see \cite{KKCh}) are being considered elsewhere.

\subsection{Preliminaries: Bilinear form and brackets}\label{FandB} A.~Lebedev proved (\cite{LeD}) that if $p=2$ and $\dim V$ is \textbf{odd, there is one equivalence class} of the \textbf{even} non-degenerate symmetric bilinear forms on the linear space $V$, whereas if $\dim V$ is \textbf{even, there are two classes}: 

\underline{type $I$)} the Gram matrix has at least one non-zero element on the main diagonal, 

\underline{type $\Pi$)} if all elements on the main diagonal of the Gram matrix are zero. 

For the normal shapes of these Gram matrices we take 
%\be
\begin{eqnarray}%\label{shapePi}
\label{sh1}\tilde \Pi_{2n}:=\begin{pmatrix}
1_2 &0\\
0& \Pi_{2n-2}\\
\end{pmatrix}\text{~~if $B$ is of type $I$ and $\dim V=2n$},\\ 
\label{sh2}\text{$\Pi_{2a}:=\begin{pmatrix}
0&1_a\\
1_a&0\\
\end{pmatrix}$ if $B$ is of type $\Pi$ and $\dim V=2a$, or}\\ 
\label{sh3}\text{$\Pi_{2a+1}:=\begin{pmatrix}
0&0&1_a\\
0&1&0\\
1_a&0&0\\
\end{pmatrix}$ if $\dim V=2a+1$ and $B$ of any type.}
\end{eqnarray}
%\ee 
We denote the Lie algebra preserving the form $B$ by $\fo_B(V)$ or $\fo_B(a)$ if $\dim V=a$. 

If $V$ is a superspace, the even form $B$ on it is the direct sum of the forms on the even and odd parts of $V$, and hence the non-degenerate symmetric forms can be of the types $B=II$, $I\Pi$, $\Pi I$, and $\Pi\Pi$ --- short for $I\oplus I$, etc. We denote the ortho-orthogonal Lie superalgebra preserving $B$ by $\fo\fo_B(V)$ or $\fo\fo_B(a|b)$ if $\sdim V=a|b$.

A.~Lebedev proved (\cite{LeD}) that for any $p$, there is one equivalence class of non-degenerate \textbf{odd} symmetric bilinear forms $B$ on $V$. In this case, $\sdim V=a|a$ and for the normal shape of $B$ one can take $\Pi_{a|a}= \Pi_{2a}$. We denote the Lie superalgebra preserving $B$ by $\fpe_B(V)$ or $\fpe_B(a)$.

$\bullet$ Define the \textit{Poisson bracket} on the space of Grassmann algebra generated by either $2n$ odd generators $\xi_1,\dots,\xi_n$ and $\eta_1, \dots,\eta_n$ (case $\Pi$), or $n$ odd generators $\theta_1,\dots,\theta_n$ (case $I$), in accordance with \eqref{sh1}--\eqref{sh3}:
\begin{eqnarray}
\label{po1}\{f,g\}_\Pi:= \sum\left(\frac{\partial f}{\partial {\xi_i}} \frac{\partial g}{\partial {\eta_i}} +\frac{\partial f}{\partial {\eta_i}} \frac{\partial g }{\partial {\xi_i}}\right);\\
\label{po2}\{f,g\}_I:= \sum\left(\frac{\partial f}{\partial {\theta_i}} \frac{\partial g }{\partial {\theta_i}}\right);\\
\label{po3}\{f,g\}_I:= \sum_{i\leq n-2}\left(\frac{\partial f}{\partial {\xi_i}} \frac{\partial g}{\partial {\eta_i}} +\frac{\partial f}{\partial {\eta_i}} \frac{\partial g }{\partial {\xi_i}}\right)+
\begin{cases}\frac{\partial f}{\partial {\theta}} \frac{\partial g }{\partial {\theta}}&\text{~for $n$ odd}\\
\sum_{i\leq 2}\left(\frac{\partial f}{\partial {\theta_i}} \frac{\partial g }{\partial {\theta_i}}\right)&\text{~for $n$ even}.
\end{cases}
\end{eqnarray}

$\bullet$ If, instead of the Grassmann algebra, we consider the algebra of truncated polynomials generated by either $2n$ even generators $q_1,\dots,q_n$ and $p_1, \dots,p_n$ (case $\Pi$), or $n$ even generators $z_1,\dots,z_n$ (case $I$), the \textit{Poisson bracket} becomes
\begin{eqnarray}
\label{po4}\{f,g\}_\Pi:= \sum\left(\frac{\partial f}{\partial {q_i}} \frac{\partial g}{\partial {p_i}} +\frac{\partial f}{\partial {p_i}} \frac{\partial g }{\partial {q_i}}\right);\\
\label{po5}\{f,g\}_I:= \sum\left(\frac{\partial f}{\partial {z_i}} \frac{\partial g }{\partial {z_i}}\right);\\
\label{po6}\{f,g\}_I:= \sum_{i\leq n-2}\left(\frac{\partial f}{\partial {q_i}} \frac{\partial g}{\partial {p_i}} +\frac{\partial f}{\partial {p_i}} \frac{\partial g }{\partial {q_i}}\right)+\begin{cases}\frac{\partial f}{\partial {z}} \frac{\partial g }{\partial {z}}&\text{~for $n$ odd}\\
\sum_{i\leq 2}\left(\frac{\partial f}{\partial {z_i}} \frac{\partial g }{\partial {z_i}}\right)&\text{~for $n$ even}.
\end{cases}\end{eqnarray}
We call the above (super)spaces with the Poisson brackets the \textit{Poisson algebras}. Observe that the space of truncated polynomials with the bracket $\{\cdot, \cdot\}_I$ is not a Lie algebra, but a Leibniz\footnote{The (left) \textit{Leibniz} algebra $L$ satisfies $[x, [y,z]]=[[x,y], z] +[y, [x,z]]$ for any $x,y,z\in L$; if, moreover, $L$ it is anti-commutative, it is a Lie algebra. Superization is immediate, via the Sign Rule.} algebra: $\{z_i, z_i\}_I=1$, not 0. The quotient modulo center is, however, a Lie algebra. For all the above Poisson brackets and the combinations thereof on $a$ even and $b$ odd indeterminates, the quotients modulo center (generated by constants) is the \textit{Lie superalgebra of Hamiltonian vector fields} $\fh_B(a|b)$, where $B=\Pi\Pi$, $\Pi I$, $I\Pi$, or $II$.

$\bullet$ We also consider the \textit{Buttin superalgebra} $\fb(a|a)$; its space is the tensor product of Grassmann algebra by the algebra of truncated polynomials generated by $n$ even generators $q_1,\dots,q_n$ and $n$ odd generators $\pi_1,\dots,\pi_n$ with the \textit{Schouten bracket} a.k.a. \textit{Buttin bracket} a.k.a. \textit{anti-bracket}:
\begin{eqnarray*}\label{Bb}
\{f,g\}_{B.b}:= \sum\left(\frac{\partial f}{\partial {q_i}} \frac{\partial g}{\partial {\pi_i}} +\frac{\partial f}{\partial {\pi_i}} \frac{\partial g }{\partial {q_i}}\right).
\end{eqnarray*}
Set $\fle(a|a)=\fb(a|a)/\fc$ and $\fh_B(a|b)=\fpo_B(a|b)/\fc$, where $\fc$ is the center (spanned by constants); $\fle(a|a)$ was introduced together with its interpretation in \cite{Le}.

On each of the above Poisson (Leibniz) algebras $\fpo_B(a|b)$ and the Buttin algebras $\fb(a|a)$, define NIS by means of the \textit{Berezin integral} = the coefficient of the highest monomial of $fg$: 
\[
B(f,g)=\int fg\vvol, \text{~~where $\vvol$ is the volume element.}
\]
 This NIS is of the same parity as the number of odd indeterminates and induces a NIS of the simple subquotients, $\fh^{(1)}_B(a|b)$ and $\fle^{(1)}(a|a)$. Clearly, $\fh^{(1)}_B(a|b)$ (resp. $\fle^{(1)}(a|a)$) have the Poisson or Leibniz algebras (resp., the Buttin algebra) as their double extensions; but they also have other double extensions.

For Lie superalgebras with Cartan matrices, see \cite{BGL}; for descriptions in terms of Cartan-Tanaka-Shchepochkina prolongations, see \cite{BGLLS1,BGLLS}; for the classification of simple Lie superalgebras, see \cite{BGL,BLLS}. 

\subsection{Preliminaries: gradings}\label{grad} A.~Lebedev argued that in the modular case, the space of roots, although a particular case of the space of weights, should be considered over $\Ree$, not over the ground field $\Kee$ as weights \textbf{are} considered in any representation except for the adjoint one, see \cite[Subsect. 4.3]{BGL}. The same applies to the notion of $\Zee$-graded Lie (super)algebras: all roots of $\fder\ \fg$, in particular, $\Zee$-gradings, should be considered over $\Ree$. (It is strange that this observation was not made ca 50 years earlier or any time later. Recently we realized that this was understood as early as in \cite{KuCh}.) This interpretation of the gradings and roots is implemented in SuperLie, see \cite{Gr}; it allowed us to accelerate computations up to 100 times: e.g., the computation time for algebras with 6 indeterminates reduces from 8 hours to 5 minutes.

\section{Background: Double extensions for $p\neq 2$ (after \cite{BKLS,ABB,Be})}

\subsection{Lemma (On a central extension)}\textup{(Lemma 3.6, page 73 in \cite{BB})} \label{L1}
\textit{Let $\fa$ be a Lie (super)algebra over a field $\Kee$, let $B_\fa$ be an $\fa$-invariant NIS on $\fa$, let $D\in\fder\ \fa$ be a derivation such that $B_\fa$ is $D$-invariant, i.e.,
\be\label{(1)}
B_\fa(Da,b)+(-1)^{p(a)p(D)}B_\fa(a,Db)=0\text{~~ for any $a,b\in \fa$. } \ee
Then,  the bilinear form 
\be\label{bilPne2}
\omega(a,b):=B_\fa(Da,b)
\ee
is a $2$-cocycle of the Lie (super)algebra $\fa$; clearly, $p(\omega)=p(B)+p(D)$.}
%\end{Lemma}

Thus, under assumptions of Lemma~\ref{L1}, we can construct a central extension $\fa_\omega$ of $\fa$ with the center spanned by an element $c$, given by cocycle $\omega$ so that $\fa_\omega/\Kee c\simeq \fa$; moreover, we can construct a semidirect sum $\fg=\fa_\omega\ltimes \Kee D$.

On the Lie (super)algebra $\fg$, define a bilinear form $B$ by setting for any $x\in\fa$
\be\label{nis}
B|_{\fa}=B_\fa, \quad B(D,c)=1, \quad B(c,x)=B(D,x)=B(c,c)=B(D,D)=0.
\ee

\subsection{Lemma (On NIS)} \label{L2} \textup{(For Lie algebras: \cite[Exercise 2.10] {K}; for Lie superalgebras: \cite[Theorem 1, page 68]{BB})}
\textit{The form $B$ defined by \eqref{nis} is a NIS.}
%\end{Lemma}

Thus constructed the Lie (super)algebra $\fg$ with NIS $B$ on it is called the \textit{double extension} or, for emphasis, \textit{D-extension}, of $\fa$. 

\subsubsection{Remark} If the derivation $D$ is inner, then $\fg$ is a direct sum of its ideal $\fa$ and a~2-dimensional commutative ideal $\fa \oplus(\Kee c\oplus \Kee D)$; this --- \textit{decomposable} --- case is not interesting.

\subsection{Lie (super)algebra $\fg$ that can be a double extension of a~ Lie (su\-per)\-algebra $\fa$}
Let $\fg$ be a Lie algebra over any field $\Kee$ or a Lie superalgebra over a field $\Kee$ of characteristic $p\ne 2$; let $B$ be a NIS on $\fg$, and $c\ne 0$ a central element of $\fg$.

The invariance
of the form $B$ implies that $B(c,[x,z])=0$ for any $x,z\in\fg$, i.e., the space $c^\perp$ contains the commutant $\fg^{(1)}=[\fg,\fg]$ of $\fg$, and hence is an ideal. Since the form $B$ is non-degenerate, the codimension of this ideal is equal to 1.

If $B(c,c)\ne 0$, then the Lie (super)algebra $\fg$ is just a direct sum $\fg=\Kee c\oplus c^\perp$. 

\sssbegin{Theorem}\emph{(For Lie algebras, see [MR, discussion after Prop. 1.2, page 555]; for Lie superalgebras, see  [BB, Corollary 1, page 81],  and [ABB1, Corollary 2.15, page 238])}\label{Pr5.1}  Let $\fg$ be a Lie algebra over any field $\Kee$ or a Lie superalgebra over a field $\Kee$ of characteristic $p\ne 2$; let $B$ be a NIS on $\fg$, and $\fz(\fg)$ the center of $\fg$. 

The Lie (super)algebra $\fg$ is a double extension of a Lie (super)algebra $\fa$ in the following cases:

\emph{(i)} $\fa$ is a Lie algebras, or Lie superalgebra with $D$ and $B$ of the same parity, and 
\[
\fg_\ev \cap \fz (\fg) \not = 0;
\]

\emph{(ii)} $\fa$ is a  Lie superalgebra with $D$ and $B$ of different parities, and
\[
\fg_\od\cap\fz(\fg)\neq 0.
\]
\end{Theorem}

\subsection{Classically restricted Lie (super)algebras (after (\cite{BLLS})} Observe that if $p=2$ there are several non-classical notions of restrictedness, see \cite{BLLS}, which we are not considering here.

\sssec{Restrictedness of Lie algebras}\label{ss-p-str} Let the ground field $\Kee$
be of characteristic $p>0$, and $\fg$ a~Lie algebra. For every $x\in
\fg$, the operator $(\ad_x)^{p}$ is a derivation of
$\fg$. If this derivation is an inner one, i.e., there is a map
(called \textit{$p$-structure}) ${}[p]:\fg\tto\fg, \ x\mapsto x^{[p]}$ such that
\begin{equation} \label{restricted-3}
[x^{[p]}, y]=(\ad_x)^{p}(y)\quad \text{~for any~}x,y\in\fg,
\end{equation} 
\begin{equation}\label{restricted-1}
(ax)^{[p]}=a^px^{[p]}\quad \text{~for any~}a\in\Kee,~x\in\fg,
\end{equation}
\begin{equation} \label{restricted-2}
(x+y)^{[p]}=x^{[p]}+y^{[p]}+\mathop{\sum}\limits_{1\leq i\leq
p-1}s_i(x, y) \quad\text{~for any~}x,y\in\fg,
\end{equation}
where $is_i(x, y)$ is the coefficient of $\lambda^{i-1}$ in
$(\ad_{\lambda x+y})^{p-1}(x)$,
then the Lie algebra
$\fg$ is said to be \emph{restricted} or \emph{having a
$p$-structure}.

\sssbegin{Remark}\label{Uniq} If the Lie algebra $\fg$ is centerless,
then the condition \eqref{restricted-3} implies  \eqref{restricted-1} and  \eqref{restricted-2}.
A $p$-structure on a given Lie algebra $\fg$ does not have to be unique; all $p$-structures on $\fg$ agree modulo center. Hence, on any simple Lie algebra, there is
not more than one $p$-structure.
\end{Remark}

\sssec{Restricted modules}\label{resModP} A $\fg$-module $M$ over a
restricted Lie algebra $\fg$, and the representation $\rho$ defining
$M$, are said to be \textit{restricted} or having a
\textit{$p$-structure} if
\begin{equation*} \label{restri-3}
\rho(x^{[p]})=(\rho(x))^{p}\quad \text{~for any~}x \in\fg.
\end{equation*}

\sssec{The $p|2p$-structure or restricted Lie
superalgebra}\label{SSp2pStr}

For a Lie superalgebra $\fg$ of characteristic $p>0$, let the Lie
algebra $\fg_\ev$ be restricted and
\begin{equation} \label{restr3}
[x^{[p]}, y]=(\ad_x)^{p}(y)\quad \text{~for any~}x\in\fg_\ev,~y\in\fg.
\end{equation}
This gives rise to the map (recall that the bracket of odd elements
is the polarization of the squaring $x\mapsto x^2$)
\begin{equation*}\label{2p}
{}[2p]:\fg_\od\to\fg_\ev, ~~~ x\mapsto(x^2)^{[p]},
\end{equation*}
satisfying the condition
\begin{equation*}\label{2p2}
{}[x^{[2p]},y]=(\ad_x)^{2p}(y)\quad\text{~for any~}x\in\fg_\od,~y\in\fg.
\end{equation*}
The pair of maps $[p]$ and $[2p]$ is called a $p$-\textit{structure}
(or, sometimes, a $p|2p$-\textit{structure}) on $\fg$, and $\fg$ is
said to be \textit{restricted}. It suffices to determine the $p|2p$-structure on any basis of $\fg$; on simple Lie superalgebras there are not more than one $p|2p$-structure.

\sssbegin{Remark}\label{ResQuo} If (\ref{restr3}) is not satisfied,
the $p$-structure on $\fg_\ev$ does not have to generate a
$p|2p$-structure on $\fg$: even if the actions of $(\ad_x)^p$ and
$\ad_{x^{[p]}}$ coincide on $\fg_\ev$, they do not have to coincide
on the whole of $\fg$. 
\end{Remark}

\section{Double extensions for $p=2$ (for proofs, see \cite{BeB1})}\label{S4}

\subsection{Quadratic and bilinear forms for $p=2$} A given map $q:V\rightarrow \Kee$, where $V$ is a $\Kee$-vector space, is called a \textit{quadratic form} if
\[
\begin{array}{c}
q(\lambda v)=\lambda^2 q(v)\text{ for any $\lambda \in \Kee$ and for any $v \in V$, and the map }\\
(u,v) \mapsto B_q(u,v):=q(u+v)-q(u)-q(v) \text{ is bilinear.}
\end{array}
\]
The form $B_q$ is called the \textit{polar form} of $q$. %Recall that quadratic forms over a field of characteristic 2 are classified by the
%so-called
%\textit{Arf} invariant, see \cite{D}. 
Recently, Lebedev classified non-degenerate bilinear forms over a perfect field $\Kee$ (i.e., such that $\Kee^2=\Kee$), see \cite{LeD}. This is a non-trivial result not related with a well-known classification of quadratic forms in any characteristic, because in characteristic 2, each of the arrows
\[
q \longleftrightarrow B_q
\]
is not necessarily onto, has a kernel, and two quadratic forms with different Arf invariants may have identical polar forms.

\subsection{Lie superalgebras for $p=2$} 
For $p\neq2$, superization of many notions of Linear Algebra is performed, as is now well-known, with the help of the Sign Rule. If $p=2$, \textbf{additional} conditions appear.
We recall basic definitions retaining the minus sign from definitions for $p\neq 2$: for clarity.
 
$\bullet$ A \textit{Lie superalgebra} in characteristic 2 is a
superspace $\fg=\fg_\ev\oplus\fg_\od$ over a field $\Bbb K$ such that the even part
$\fg_\ev$ is a Lie algebra, the odd part $\fg_\od$ is a
$\fg_\ev$-module made two-sided by
\textit{anti}-symmetry, and on the odd part $\fg_\od$ a \textit{squaring} is
defined as a map
\begin{equation}\label{squaring}
\begin{array}{l}
s_\fg:\fg_\od \rightarrow \fg_\ev\text{~~given by $f\mapsto f^2$ such that}\\
\text{$\bullet$  $(\lambda f)^2=\lambda^2 f^2$ for any $f\in
\fg_\od$ and $\lambda \in \Kee$, and}\\
\text{$\bullet$   the map} {}[f,g]:= (f+g)^2- f^2- g^2 \text{~~for any $f,g\in\fg_\od$}\\
%\text{(which is the bracket $\fg_\od\times \fg_\od \rightarrow \fg_\ev$) 
\phantom{...}\text{~~is a bilinear form on $\fg_\od$ with values
in $\fg_\ev$.}
\end{array}
\end{equation}
The bracket on $\fg_\ev$, as well as the action of $\fg_\ev$ on $\fg_\od$, is denoted also by the same symbol $[\cdot ,\cdot]$, or $[\cdot ,\cdot]_\fg$ for clarity. 
The Jacobi identities for three even elements, and involving one or two odd elements, are the same as for $p\neq 2$; the following new identity replaces the one with three odd elements: it involves the squaring:
\begin{equation}\label{JIS}
{}[f^2,g]=[f,[f,g]]\;\text{ for any $f\in \fg_\od$ and $g\in\fg$}.
\end{equation}

$\bullet$ Desuperizing (retain only the bracket) we get a $\Zee/2$-graded Lie algebra $\fg$.

$\bullet$ For any Lie superalgebra $\fg$ in characteristic 2, its \textit{derived algebras} are
\be\label{prim}
\fg^{(0)}: =\fg, \quad
\fg^{(i+1)}=[\fg^{(i)},\fg^{(i)}]+\Kee\{f^2\mid f\in (\fg^{(i)})_\od\}.
\ee
A linear map $D:\fg\rightarrow \fg$ is called a \textit{derivation} of the Lie superalgebra $\fg$ if \begin{eqnarray}
\label{Der1} D([f,g])&=&[D(f),g]+[f,D(g)]\quad \text{for any $f,g\in \fg$}
\end{eqnarray}
and, \textbf{additionally} (generalization of \eqref{JIS}, where $D=\ad_g$), 
\begin{eqnarray}\label{Der2} D(f^2)&=&[D(f),f]\quad \text{for any $f\in \fg_\od$}.
\end{eqnarray}

%$\bullet$ It is worth noticing that condition (\ref{Der2}) implies condition (\ref{Der1}) if $f,g\in \fg_\od$.
%In showing that $D$ is a derivation, we will not check condition (\ref{Der1}) for every $f$ and $g$ odd as long as condition $(\ref{Der2})$ has been checked.

We denote the space of all derivations of $\fg$
% $(\fg,[\cdot ,\cdot]_\fg,s_\fg)$
by $\fder (\fg)$.

An even linear map $\rho: \fg\tto\fgl(V)$ is a \textit{representation
of the Lie superalgebra} $\fg$ in the superspace $V$, called \textit{$\fg$-module}, if
\begin{equation}\label{repres}
\begin{array}{l}
\rho([f, g])=[\rho (f), \rho(g)]\text{ for any $f, g\in
\fg$};\\ 
\text{ and, \textbf{additionally}, }
\rho (f^2)=(\rho (f))^2\text{~for any $f\in\fg_\od$.}
\end{array}
\end{equation}
We say that a bilinear form 
 $B$ on $\fg$ is \textit{symmetric} if 
\[
\begin{array}{l}
\text{$B(f,g)=B(g,f)$ for any $f,g\in\fg$; and, }\\
\text{\textbf{additionally}, $B(f,f)=0$ for any $f\in \fg_\od$}.
\end{array}
\]
We say that $D$ \textit{preserves} $B$ if
\be\label{extraB}
\begin{array}{l}
\text{$B(D(f),g)+B(f,D(g))=0$ for any $f,g\in \fg$; and, }\\
\text{\textbf{additionally}, $B(D(f), f)=0$ for any $f\in \fg_\od$.}
\end{array}
\ee

We denote the \textit{NIS-superalgebra} with a NIS $B_\fg$ by $(\fg, B_\fg)$.
A NIS-superalgebra $(\fg, B_\fg)$ is said to be \textit{decomposable} if it can be decomposed into direct sums of ideals, namely $\fg=\oplus I_i$, such that all $I_i$ are orthogonal to each other.

%%%%%%%%%%%%%%%%%%
\subsection{The case where $B_\fa$ is even} Let $\fg:=\mathscr{K} \oplus \fa \oplus \mathscr{K} ^*$ as spaces. %, $\mathscr{K} :=\Kee c$ and $\mathscr{K}^* :=\Kee D$.

 \subsubsection{$D_\ev$-extensions. The quadratic form $q$} {}~{}
 
\sssbegin{Theorem}\label{MainTh} Let $(\fa, B_\fa)$ be a NIS-superalgebra, $p(B_\fa)=\ev$, let $D\in \fder_\ev(\fa)$ preserve $B_\fa$.
Let $q: \fa_\od \rightarrow \Bbb K$ be a quadratic form; let its polar form $B_q$ satisfy
\begin{eqnarray}
\label{D3}
B_q(a,b)&=&B_\fa(a, D(b)) \text{ for any $a,b \in \fa_\od$}.
\end{eqnarray}

Then,  there exists a NIS-superalgebra structure on $\fg$, defined as follows, cf. Subsection~$\ref{L1}$ and \eqref{nis}. The squaring is given by
\[
s_\fg(a):= s_\fa(a)+q (a) c \qquad \text{ for any } a\in \fg_\od \; (=\fa_\od).
\]
The bracket on $\fg$ is defined for any $a,b\in \fa$ as follows:
\be\label{*}
 [a,b]_\fg:=[a,b]_\fa+B_\fa(D(a),b)c, \quad [D,a]_\fg:=D(a)\text{ and $[c,\fg]_\fg=0$}.
\ee
The NIS $B$ on $\fg$ is given by formulas \eqref{nis} with one modification: $B(D, D)$ can be arbitrary.
%\end{eqnarray*}
%Moreover, the form $B_\fg$, obviously even, is invariant on $\fg$.
\end{Theorem}

What is the difference of this Theorem from its analog for $p\ne 2$? 

\underline{First}, the additional condition in eq.~ \eqref{extraB}: $B(D(f),f)=0$. For $p\ne 2$, it follows from eq.~ \eqref{extraB}, whereas for $p=2$ it should be required as a part of the invariance condition. 

\underline{Second}, the reasons for the quadratic form $q$ to appear. Indeed, for $p\ne 2$, to determine a central extension, we only need a cocycle, whereas for $p=2$, we have to adjust the squaring to match this cocycle and this is precisely what the form $q$ is needed for. 

\textbf{Everything else is the same as for $p\ne 2$, same as with the desuperisations of $\fa$}, when only even derivations and forms remain, i.e., no extra conditions.

We call the Lie superalgebra $(\fg, B)$ constructed in Theorem \ref{MainTh} a \textit{$D_\ev$-extension} of $(\fa, B_{\fa})$ by means of $D$ and $q$.

Now, let the \textit{\lq\lq special center" of $\fg$ relative to $B_\fg$} be
\begin{equation}\label{speCen}
\fz_s(\fg):=\fz(\fg)\cap s_\fg(\fg_\od)^{\perp}, \text{~~where $s_\fg(\fg_\od):=\Span\{x^2\mid x\in\fg_\od\}$}.
\end{equation}
Observe that $\fz_s(\fg)_\od=\fz(\fg)_\od$ and $\fz_s(\fg)_\ev=\fz(\fg)_\ev \cap s_\fg(\fg_\od)^{\perp}$. Moreover, $\fz_s(\fg)$ is not necessarily an ideal.

\sssbegin{Proposition}
\label{Rec1}
Let $(\fg,B_\fg)$ be an indecomposable NIS-superalgebra. Let $\fz_s(\fg)_\ev\not=\{0\}.$ Then,  $(\fg,B_\fg)$ is obtained as an $D_\ev$-extension from a NIS-superalgebra $(\fa,B_\fa)$ with $p(B_\fa)=\ev$. 
\end{Proposition}

What is the difference of the case where $p=2$ from that where $p\ne 2$? The ``indecomposability'' condition can be weakened to the existence of a nonzero even element $c$ belonging to the intersection  $[\fg,\fg]\cap s_\fg(\fg_\od)^\perp$. Now we see that the only extra condition is related to squaring and is precisely due to the fact that $c^\perp$ must be an ideal.

\subsubsection{$D_\od$-extensions. The element $A\in \fa_\ev$}

\sssbegin{Theorem}\label{MainTh2} Let $(\fa,B_\fa)$ be a NIS-superalgebra, $p(B_\fa)=\ev$. Let $D\in \fder_\od(\fa)$ preserve $B_\fa$; let $A\in \fa_\ev$ satisfy the following conditions \emph{(see \cite{ABB,Be})}:
\begin{eqnarray}
\label{2D2} D^2=\ad_{A};\ \ D(A)=0.
\end{eqnarray}
Then,  there exists a NIS-superalgebra structure on $\fg$, defined as follows. The squaring is given by
\[
s_\fg(rc+a+t D):= s_\fa(a) +t D(a) + t^2 A\qquad \text{ for any $a\in \fa_\od$ and $r,t\in\Kee$}.
\]
The bracket is given by eq.~\eqref{*}. The even NIS $B$ on $\fg$ is given by formulas \eqref{nis}. Besides, in this case we only get $D_\od$-extensions.
\end{Theorem}

We call the NIS-superalgebra $(\fg,B)$ constructed in Theorem
\ref{MainTh2} a \textit{$D_\od$-extension} of $(\fa, B_{\fa})$ by means of $D$ and $A$.

\parbegin{Proposition}
\label{Rec2}
Let $(\fg, B)$ be a NIS-superalgebra, $p(B)=\ev$. Define the cone
\[
\mathscr{C}(\fg, B):=\{x\in \fg_\od\; | \; B(s_\fg(x), s_\fg(y))=0 \ \ \text{ for any } y\in \fg_\od \}.
\]
Let $\fz(\fg)_\od \cap \mathscr{C} (\fg, B)\not=\{0\}$. Then,  $(\fg, B)$ is obtained from a NIS-superalgebra $(\fa,B_\fa)$ as  a $D_\od$-extension. % of dimension $\dim(\fg) - 2$.
\end{Proposition}

\subsection{The case where $B_\fa$ is odd} Here $\fg:=\mathscr{K} \oplus \fa \oplus \mathscr{K} ^*$ as spaces, where $\mathscr{K} :=\Kee c$ and $\mathscr{K}^* :=\Kee D$.

\subsubsection{$D_\ev$-extensions}

\sssbegin{Theorem}\label{MainTh4} Let $(\fa, B_\fa)$ be a NIS-superalgebra, $p(B_\fa)=\od$. Let $D\in \fder_\ev(\fa)$ preserve $B_\fa$. Then,  there exists a NIS-superalgebra structure on $\fg$, 
defined as follows. The squaring is given by
\[
s_\fg(a+\mu c):= s_\fa(a) \qquad \text{ for any $a\in \fa_\od$ and $\mu\in\Kee$}.
\]
The bracket is given by eq.~\eqref{*}. The NIS $B$ on $\fg$ is given by formulas \eqref{nis}.
\end{Theorem}

We call the NIS-superalgebra $(\fg,B)$ constructed in Theorem \ref{MainTh4} the \textit{$D_\ev$-extension} of $(\fa, B_{\fa})$ by means of $D$.

\parbegin{Proposition}\label{Rec4}
Let $(\fg, B)$ be an indecomposable NIS-superalgebra, $p(B)=\od$. 
If $s_\fg(\fz(\fg)_\od) \cap s_\fg(\fg_\od)^\perp\not=\{0\}$, then $(\fg, B)$ is obtained from a NIS-superalgebra $(\fa,B_\fa)$ as a $D_\ev$-extension. 
\end{Proposition}

\subsubsection{$D_\od$-extensions. The quadratic form $q$, and the elements $A\in\fa_\ev$ and $m\in\Kee$}

\parbegin{Theorem}\label{MainTh3} Let $(\fa, B_\fa)$ be a NIS-superalgebra, $p(B_\fa)=\od$. Let also $A\in \fa_\ev$ and $D\in \fder_\od(\fa)$ preserving $B_\fa$ satisfy the conditions \eqref{2D2}.

Let $q$ be a quadratic form on $\fa_\od$ such that $B_q(a,b)=B_\fa(D(a),b)$. 

Then,  there exists a NIS-superalgebra structure on $\fg$, 
defined as follows. The squaring is given by
\begin{eqnarray*}
s_\fg(a+\mu D):= s_\fa(a) +(\mu^2m+q(a))c+ \mu^2 A+\mu D(a) \\ 
\text{for any $a\in \fa_\od$ and $\mu\in\Kee$, and some $m\in\Kee$}.
\end{eqnarray*}
The bracket is given by eq.~\eqref{*}. The NIS $B$ on $\fg$ is given by formulas \eqref{nis}.
\end{Theorem}

We call the NIS-superalgebra $(\fg, B)$ constructed in Theorem \ref{MainTh3} a \textit{$D_\od$-extension} of $(\fa, B_{\fa})$ by means of $D$, $q$, $A$, and $m$.

\parbegin{Proposition}
\label{Rec3}
Let $(\fg, B)$ be an indecomposable NIS-superalgebra, $p(B)=\od$. If $\fz(\fg)_\ev\not=\{0\}$, then $(\fg, B)$ is obtained as a $D_\od$-extension from a NIS-superalgebra. % $\fa$. 
\end{Proposition}

\subsection{Summary: What is needed to construct a double extension of the NIS-(super)algebra $\fa$} We need $q,A,m$ according to Table \eqref{Summ}.
\begin{equation}\label{Summ}
\begin{tabular}{|c|c|c|}\hline
\diagbox{$D$}{$B_\fa$}&even&odd\\ \hline
even& $q$ & -- \\\hline
odd & $A$ & $q,A,m$\\\hline
\end{tabular}
\end{equation}

%%%%%%%%%%%%%%%%%%%%%%%%%%%%%%%%%%%%%%%%%%%
\section{Isomorphisms, and equivalence classes of derivations (after \cite{BeB1})}\label{Siso}
%%%%%%%%%%%%%%%%%%%%%%%%%%%%%%%%%%%%%%%%%%%
For a NIS-superalgebra $\fa$ with NIS $B_\fa$, denote by $(\fg, B_\fg)$ (resp. $(\tilde \fg, \tilde B_\fg)$) the double extension of $\fa$ by means of a derivation $D$ (resp. $\tilde D$), i.e., $\fg:=\mathscr{K} \oplus \fa \oplus \mathscr{K}^*$, where $\mathscr{K}=\Kee c$ and $\mathscr{K}^*=\Kee D$ (resp. $\tilde \fg:=\tilde {\mathscr{K}} \oplus\fa \oplus \tilde {\mathscr{K}^*}$, where $\tilde {\mathscr{K}}=\Kee \tilde c$ and $\tilde {\mathscr{K}^*}=\Kee \tilde D$). 

\textbf{If $p=2$}, we \textbf{additionally} require a quadratic form $q$ (resp. $\tilde q$) if the center $c$ is even, while in the case of $D_\od$-extensions, we need, moreover, an element $A\in\fa_\ev$, and sometimes $m\in\Kee$ (resp. $\tilde A$ and $\tilde m$) described in the previous section. 

An \textit{isomorphism} of NIS-superalgebras $\pi: \fg\tto\tilde \fg$ is a NIS-preserving isomorphism of Lie superalgebra structures: 
\[
\begin{array}{l}
\pi([f,g]_\fg)=[\pi(f), \pi(g)]_{\tilde \fg},\ \ \
B_{\tilde \fg}(\pi(f), \pi(g))=B_\fg(f,g)\text{~~for any $f,g\in \fg$};\\
\text{if $p=2$, we additionally require $ \pi(s_\fg(f)) = s_{\tilde \fg} (\pi (f))$ for any $f\in \fg_\od$}.
\end{array}
\]

\subsection{The case where $B_\fa$ is even}

\sssbegin{Theorem} \label{Isom1} Let $p(D)=p(\tilde D)=\ev$.
Let $\pi_0$ be an automorphism of $(\fa, B_\fa)$. Let $\lambda\in\Kee^{\times}$ and let $y\in \fa_\ev$ satisfy the following conditions:
\[\label{Ca}
\begin{array}{l}
\tilde q(a)=\lambda q \circ \pi_0^{-1}(a)+B_\fa(y, s_\fa \circ \pi_0^{-1}(a)) \quad \text{for any $a\in\fa_\od$}; \\
\pi_0^{-1}\tilde D\pi_0\vert_{\fa_\ev}=(\lambda D+\ad_{y})\vert_{\fa_\ev};\\ 
B_\fg(D, D)=\lambda^{-2}( B_\fa(y, y)+B_\fg(\tilde D, \tilde D)).
\end{array}
\]
Then,  there exists an isomorphism $\pi: \fg\tto \tilde{\fg}$ given by the formulas
\[
\begin{array}{lcl}
\pi(a)&=& \pi_0(a)+ B_\fa(y,a) \tilde c \quad \text{for any $a\in\fa$;}\\
\pi(c)&=&\lambda \tilde c;\\
\pi(D)&=&\lambda^{-1}(\tilde D+\pi_0(y))+\nu \tilde c,\; \text{where $\nu\in\Kee$ is arbitrary.}\\
\end{array}
\]

If $[D]=[\tilde D]$ in $\mathrm{H}^1_\ev(\fa; \fa)$, then they define isomorphic double extensions of $\fa$.
\end{Theorem}

\sssbegin{Theorem} \label{Isom2} Let $p(D)=p(\tilde D)=\od$.
Let $\pi_0$ be an automorphism of $(\fa, B_\fa)$. Let $\lambda\in\Kee^{\times}$ and $y\in\fa_\od$ satisfy the following condition:
\[
\begin{array}{lcl}
\pi_0^{-1}\tilde D\pi_0\vert_{\fa}&=&(\lambda D+\ad_{y})\vert_{\fa}.
\end{array}
\]
 Then,  there exists an isomorphism $\pi: \fg\tto \tilde{\fg}$ 
 given by the formulas
\[
\begin{array}{lcl}
\pi(a)&=& \pi_0(a)+ B_\fa(y,a) \tilde c \quad \text{for any $a\in\fa_\od$};\\
\pi\vert_{\fa_\ev}&=&\pi_0\vert_{\fa_\ev};\\
\pi(c)&=&\lambda \tilde c;\\
\pi(D)&=&\lambda^{-1}(\tilde D+\pi_0(y))+\nu \tilde c, \; \text{where $\nu$ is arbitrary};\\
\tilde A&=& \lambda^2 \pi_0(A)+s_\fa(\pi_0(y))+ \lambda \pi_0(D(y)).
\end{array}
\]
 
If $[D]=[\tilde D]$ in $\mathrm{H}^1_\od(\fa; \fa)$, then they define isomorphic double extensions of $\fa$.
\end{Theorem}

\subsection{The case where $B_\fa$ is odd}

\sssbegin{Theorem} \label{Isom3}
Let $p(D)=p(\tilde D)=\od$. Let $\pi_0$ be an automorphism of $(\fa, B_\fa)$. Let $\lambda\in\Kee^{\times}$ and $y\in \fa_\od$ satisfy the following conditions:
\begin{eqnarray*}
\tilde q(a)=\lambda q \circ \pi_0^{-1}(a)+B_\fa(y, s_\fa \circ \pi_0^{-1}(a)) \quad \text{for any $a\in\fa_\od$}; \\
\pi_0^{-1}\tilde D\pi_0(a)=\lambda D(a)+\ad_{y}(a) \quad \text{for any $a\in\fa_\ev$}; \\
\tilde A=\lambda^2 \pi_0(A)+\lambda \pi_0(D(y))+s_\fa(\pi_0(y));\\
\tilde m=\lambda^2 q(y)+\lambda B_\fa(y,s_\fa(y)+\lambda A)+\lambda^{3}m.
\end{eqnarray*}
Then,  there exists an isomorphism $\pi: \fg\tto \tilde{\fg}$ given by the formulas
\be\label{umn}
\begin{array}{lcl}
\pi(a)&=& \pi_0(a)+ B_\fa(y,a) \tilde c \quad \text{for any $a\in\fa$};\\
\pi(c)&=&\lambda \tilde c;\\
\pi(D)&=&\lambda^{-1}(\tilde D+\pi_0(y)).\\
\end{array}
\ee
 
If $[D]=[\tilde D]$ in $\mathrm{H}^1_\od(\fa; \fa)$, then they define isomorphic double extensions of $\fa$.
\end{Theorem}

\sssbegin{Theorem} \label{Isom4}
Let $p(D)=p(\tilde D)=\ev$. Let $\pi_0$ be an automorphism of $(\fa, B_\fa)$. Let $\lambda\in\Kee^{\times}$ and $y\in \fa_\ev$ satisfy the following conditions:
\begin{eqnarray*}\label{4Cd}
\pi_0^{-1}\tilde D\pi_0(a)&=&\lambda D(a)+\ad_{y}(a) \quad \text{for any $a\in\fa$}.
\end{eqnarray*}
Then,  there exists an isomorphism $\pi: \fg\tto \tilde{\fg}$ given by the formulas
\eqref{umn}.
If $[D]=[\tilde D]$ in $\mathrm{H}^1_\ev(\fa; \fa)$, then they define isomorphic double extensions of $\fa$.
\end{Theorem}

%%%%%%%%%%%%%%%%%%%%%%%%%%%%%%%%%%%%%%%%%%%%%

\section{The exceptional cases: $\fa=\mathfrak{h}_B^{(1)}(a | b)$ for $a+b=4$ and 
$\fa=\mathfrak{le}^{(1)}(2 |2)$}\label{SEx}

Recall that the $\Zee$-grading of a given vectorial Lie superalgebra is called \textit{standard} if the degree of each indeterminate $X_i$ is equal to 1. This \textit{grading} is given by the Euler operator $\sum X_i\del_{X_i}$. The \textit{roots} are given with respect to the maximal torus, spanned by $q_ip_i$, $\eta_i\xi_i$, and $q_i\pi_i$ in the cases
$\mathfrak{po}_\Pi(a|0)$, $\mathfrak{po}_\Pi(0|a)$, and $\mathfrak{le}(a|a)$, respectively. According to Subsection~\ref{grad} the degree and roots are considered over $\Ree$, assuming that the $i$th coordinate of the weight of $\xi_i$, $\pi_i$, and $p_i$ is equal to 1, whereas the $i$th coordinate of the weight of $\eta_i$, and $q_i$ is equal to $-1$, other coordinates being 0.

The standard $\Zee$-grading of $\fa=\mathfrak{h}_B^{(1)}(a | b)$ (resp.
$\mathfrak{le}^{(1)}(a |a)$) is symmetric, i.e.,
\be\label{symm}
\fa=\mathop{\oplus}\limits_{-1\leq i\leq 1} \fa_i\text{~~with $\fa_{-1}\simeq\fa_1$ (resp. $\fa_{-1}\simeq\Pi(\fa_1)^*$) as $\fa_0$-modules},
\ee
only if $a+b=4$ (resp. $a=2$). That is why some of the double extensions in these cases have no analogs in the generic cases: they appear due to symmetry.

Observe the following isomorphisms, see \cite{KrLe}, that explain a result of \cite{BBH}:
\[
\fpe^{(2)}(3)\simeq\fh_\Pi^{(1)}(0|4)\simeq\fpsl(2|2)\simeq\fh_{\Pi\Pi}^{(1)}(2;\One |2)\simeq\fsvect^{(1)}(1;\One|2).
\]

\subsection{Outer derivations} Let the superscript of the derivation $D$ be its degree, let the subscript be its weight or a monomial ($x$ or $b$) or a label $\theta$.

\subsubsection{$\mathfrak{h}_{\Pi}^{(1)}(0|4)$}\label{sss04} One outer derivation in each of degrees $\pm2$ and 5 derivations of degree 0, whose weights are $(\pm 2, 0)$, $(0, \pm 2)$ and $(0, 0)$, see \cite{BeB1}.

\subsubsection{$\mathfrak{h}_{\Pi}^{(1)}(4 |0)$}\label{sss40} One outer derivation in each of degrees $\pm2$ and 5 derivations of degree 0, whose weights are $(\pm 2, 0)$, $(0, \pm 2)$ and $(0, 0)$. This case is the desuperization of the case $\mathfrak{h}_{\Pi}^{(1)}(0|4)$; the shapes of cocycles are identical  to those of $\mathfrak{h}_{\Pi}^{(1)}(0|4)$.

\subsubsection{$\mathfrak{h}_{\Pi\Pi}^{(1)}(2 |2)$}\label{sss22} One outer derivation in each of degrees $\pm2$ and 5 derivations of degree 0, whose weights are $(\pm 2, 0)$, $(0, \pm 2)$ and $(0, 0)$. This case is a partial desuperization of the case $\mathfrak{h}_{\Pi}^{(1)}(0|4)$; the shapes of cocycles are identical to those of $\mathfrak{h}_{\Pi}^{(1)}(0|4)$.

\subsubsection{$\mathfrak{h}_{II}^{(1)}(2 |2)$, $\mathfrak{h}^{(1)}(1 |3)$, $\mathfrak{h}^{(1)}(3 |1)$, $\mathfrak{h}_{I\Pi}^{(1)}(2 |2)$, $\mathfrak{h}_{\Pi I}^{(1)}(2 |2)$}\label{sssAB} One outer derivation in each of degrees $\pm2$.
The apparently missing Euler operator is in these cases an inner derivation. The same reason causes absence of certain derivations in cases $\mathfrak{h}_{I}^{(1)}(0|4)$, and $\mathfrak{h}_{I}^{(1)}(4 |0)$, and in these cases as compared with the case $\mathfrak{h}_{\Pi}^{(1)}(0|4)$.

\subsubsection{$\mathfrak{h}_{I}^{(1)}(0|4)$}\label{sss04I} One outer derivation in each of degrees $\pm2$ and 4 derivations of degree 0, whose weights are $\pm 2$ and two derivations of weight $0$.

\subsubsection{$\mathfrak{h}_{I}^{(1)}(4 |0)$}\label{sss40Ii} This is the desuperization of $\mathfrak{h}_{I}^{(1)}(0|4)$; derivations are of the same shape.

\subsubsection{$\mathfrak{le}^{(1)}(2 |2)$}\label{sssLE} The shape of derivations is identical to those of $\mathfrak{h}_{\Pi}^{(1)}(0|4)$, but the $D_{b}^{(0)}$ are odd.

\subsection{The double extensions} \label{DEevB} Note: in all exceptional cases $p(B_\fa)=\ev$.

\textbf{Case $\fh$}.\label{DEevB4} For the proof for the cocycle $D^{(-2)}$ in case $\fa=\fh_\Pi^{(1)}(0|4)$, see \cite{BeB1}; the result for $D^{2}$ is isomorphic due to the symmetry \eqref{symm}. The idea of the proof is identical for the other derivations. The condition $B_\fa^TD^{(2)}=D^{(2)}B_\fa$ is easily seen. Besides, since $D^{(2)}$ acts by zero on $\fa_\ev$, it follows that
\[
B_{\fa}(D^{(2)}(f),f)=0\text{~~ for any $f \in \fa_\ev$. }
\]

Table  \eqref{h04i} (resp. \eqref{h04}) show existence of the quadratic form $q$ (resp. element $A$) associated with each derivation $D$, and names of the respective double extensions DE. \begin{equation}\label{h04i}
\footnotesize
\renewcommand{\arraystretch}{1.4}
\begin{tabular}{|c|c|c|c|} \hline
Derivation & $q$&
DE of $\fh_I^{(1)}(0|4)$&
DE of $\fh_I^{(1)}(4|0)$\\
\hline
$D_{0}^{(0)}$ & yes & $\widehat{\mathfrak{po}}_I(0|4)$& $\widehat{\mathfrak{po}}_I(4|0)$\\ \hline
$D_b^{(0)}$ & yes & $\widetilde{\mathfrak{po}}_I(0|4)$& $\widetilde{\mathfrak{po}}_I(4|0)$\\ \hline
$D_\theta^{(0)}$ & $-$ & $-$& $-$\\ \hline
$D^{(\pm 2)}$ & yes & $\mathfrak{po}_I(0|4)$& $\mathfrak{po}_I(4|0)$\\ \hline
\end{tabular}
\end{equation}

In the 6th column of \eqref{h04} $B=II$, or $I\Pi$, or $\Pi I$. 
\begin{equation}\label{h04}
\footnotesize
\renewcommand{\arraystretch}{1.4}
\begin{tabular}{|c|c|c|c|c|c||c|c|} \hline
$D$ & $q$&
$\fh_\Pi^{(1)}(0|4)$&
$\fh_\Pi^{(1)}(4|0)$&
$\fh_{\Pi\Pi}^{(1)}(2|2)$&
$\fh_B^{(1)}(2|2)$& $A$&
$\fle^{(1)}(2|2)$\\
\hline
$D_{b}^{(0)}$ & yes & $\widetilde{\mathfrak{po}}_\Pi(0|4)$& $\widetilde{\mathfrak{po}}_\Pi(4|0)$& $\widetilde{\mathfrak{po}}_{\Pi\Pi}(2|2)$& $-$ &$0$ &$\tilde\fb(2|2)$\\ 
\hline
$D_0^{(0)}$ & yes & $\fgl(2|2)$& $\fgl(4|0)$&$\fgl(2|2)$& $-$& $-$&$-$\\ \hline
$D^{(2)}$ & yes & $\mathfrak{po}_\Pi(0|4)$& $\mathfrak{po}_\Pi(4|0)$& $\mathfrak{po}_{\Pi\Pi}(2|2)$& $\mathfrak{po}_B(2|2)$& $0$ &$\fb(2|2)$\\ \hline
\end{tabular}
\end{equation}

\textbf{Case $\fle$}.\label{le} Compare with the general case in Subsection~\ref{n_n_Le}.

\sssbegin{Lemma} \emph{1a)} The double extensions $\mathfrak{gl}(2|2)$, $\mathfrak{po}_\Pi(0|4)$,  and $\widetilde{\mathfrak{po}}_\Pi(0|4)$ are not isomorphic to each other.

\emph{1b)} The double extensions $\mathfrak{gl}(2|2)$, $\mathfrak{po}_{\Pi\Pi}(2|2)$,  and $\widetilde{\mathfrak{po}}_{\Pi\Pi}(2|2)$ are not isomorphic to each other.

 \emph{1c)} The double extensions $\mathfrak{gl}(4)$, $\mathfrak{po}_\Pi(4|0)$,  and $\widetilde{\mathfrak{po}}_\Pi(4|0)$ are not isomorphic to each other.

\emph{2)}  $\mathfrak{po}_I(0|4)$, $\widetilde{\mathfrak{po}}_I(0|4)$ and $\widehat{\mathfrak{po}}_I(0|4)$ are not isomorphic to each other.

\emph{3)}  $\tilde{\mathfrak{b}}(2|2) \not\simeq \mathfrak{b}(2|2)$.
\end{Lemma}

\begin{proof} Let the super-rank of the operator $A$ in the superspace $V$ be $\sdim(V/\Ker A)$. The claims are true for the following reasons.

1) Because the super-rank of $D_b^{(0)}$ is $(2,2)$ and no element in  $\mathfrak{po}$ and $\mathfrak{gl}$ has such rank, it follows that $\widetilde{\mathfrak{po}}$ is not isomorphic to the other two double extensions which are manifestly non-isomorphic to each other.

2) Because $\rk D_0^{(0)}=8$ in $\widehat{\mathfrak{po}}_I(0|4)$ and no element in $\mathfrak{po}_I(0|4)$ and in $\widetilde{\mathfrak{po}}_I(0|4)$ has such rank, it follows that $\widehat{\mathfrak{po}}_I(0|4)$ is not isomorphic to the other two double extensions. $\mathfrak{po}_I(0|4)$ and $\widetilde{\mathfrak{po}}_I(0|4)$ are non-isomorphic for the same reason as in 1).

3) Because $\rk \ad_{p_i}= \rk \ad_{\pi_i} = 7$ in $\tilde{\mathfrak{b}}(2|2)$ and no element in $\mathfrak{b}(2|2)$ has such rank. \end{proof} 

\section{The general cases. Outer derivations}\label{SgenOUT}

Let $\bar X$ be the product of all indeterminates. Derivations of $\fa$ will be called \textit{equivalent} if 
they lie on one orbit
of the group of automorphisms of $\fder\ \fa$. The derivations of $\fa$ are 1-cocycles of $\fa$ that is why we used the notation $\hat x$ instead of $x^*$ speaking about the element of $C^1(\fa)$ --- 1-cochain of $\fa$ with trivial coefficients: in case we have to multiply such cochains, we should recall that $p(\hat x)=p(x)+\od$.

\subsection{$\mathfrak{out}(\fh_\Pi^{(1)}(0|2n))$}\label{0_2nPi}
Notation for convenience: \begin{enumerate}
\item tags ``($\fle$: always odd)'' and  ``($\fle$: odd with $n$)'' describing the parity of the cocycle are used in Subsection \ref{n_n_Le}:
	\item Basis: $B = \{\eta_1,\dots ,\eta_n; \xi_1,\dots ,\xi_n\}$.
	\item The set of products of any $k$ basis elements (Choose):
	
	 $Ch_B(k)=\{x_1\dots x_k \mid x_i\in B\}$; 
	 
	 $Ch_\xi(k)=\{x_1\dots x_k \mid x_i\in \{\xi_1,\dots ,\xi_n\}\}$; 
	 
	 $Ch_\eta(k)=\{x_1\dots x_k \mid x_i\in \{\eta_1,\dots ,\eta_n\}\}$. For example, 
\[
Ch_B(2) = \{\eta_i\eta_j, \xi_i\xi_j \mid\text{for all~~}i\ne j\le n \} \cup \{\eta_j\xi_i\mid i,j\le n\}.
\]
	\item Index of elements: $Ind(x) = \text{The set of indices of $x$ in terms of ~} \eta, \xi$. 
	
For example, $Ind(\eta_1\xi_3) = \{1,3\}$, $Ind(\eta_2\xi_5\xi_7) = \{2,5,7\}$.
	\item The switch symbol: $S(\xi_i) = \eta_i$ and $S(\eta_i) = \xi_i$.
	\item Let $C(x)$ denote the set of all monomials each of which is a multiple of $x$; 
	\item let $\tilde{C}(x)$ denote the set of all monomials each not a multiple of $x$.
	\item Let $O$ denote the set of all monomials of odd degree.
\end{enumerate}
% Cocycles (tags ``$\fle$: odd'' are used in Subsections \ref{le2n}):
% % Deg$=-2$: 
% 	\begin{equation}
% 		D^{(-2)}=\sum_{1\leq i\leq2n}\ \sum_{x\in Ch(i)}\ \sum_{j\notin {Ind(x)}} x\otimes (\widehat{x\eta_j\xi_j}).\text{~~($\fle$: odd)}
% 	\end{equation}

% Deg=0:
	
	(a) Deg=0. There are $2n$ equivalent derivations: for any $b\in B$, we have
		\begin{equation}\label{b}
		D_b^{(0)}:=\sum_{0\leq i\leq 2n-2}\ \sum_{x\in Ch_{B\setminus\{b\}}(i)} (bx)\otimes(\widehat{S(b)x}).\text{~~($\fle$: always odd)}
		\end{equation}
		
		(b) Deg=0. One particular derivation --- the Euler operator
		\begin{equation}
		D_{0}^{(0)}=\sum_{1\leq i\leq n-1} \ \sum_{0\leq j\leq 2i+1}\ \sum_{x\in Ch_\xi(2i)}\ \sum_{y \in Ch_\eta(j)} (xy)\otimes(\widehat{xy}).
		\end{equation}
	%\end{enumerate}

	(c)Deg$=2n-2$.
	\begin{equation}\label{2n-2}
		D^{(2n-2)}:=\sum_{x\in B} \frac{\partial \bar X}{\partial x} \otimes(\widehat{S(x)}). \text{~~($\fle$: odd with $n$)}
	\end{equation}
%\end{enumerate}

\subsection{$\mathfrak{out}(\fh_\Pi^{(1)}(2n|0))$}\label{2n_0Pi}The desuperization of the case $\mathfrak{out}(\fh_\Pi^{(1)}(0|2n))$; the same cocycles.

%\subsection{$\mathfrak{out}(\fh_{\Pi}^{(1)}(0|2b))$}\label{0_2bI}

\subsection{$\mathfrak{out}(\fh_{\Pi\Pi}^{(1)}(2a|2b))$}\label{2a_2bPiPi} A partial desuperization of $\mathfrak{out}(\fh_\Pi^{(1)}(0|2a+2b))$; the same cocycles.

\subsection{$\mathfrak{out}(\fh_{I}^{(1)}(0|2n))$}\label{0_2b_I} Cocycles:
	
	(a) Deg=0. There are $2n$ equivalent derivations: for any $b\in B$, we have
		\begin{equation}
		D_b^{(0)}:=\sum_{0\leq i\leq 2n-2}\ \sum_{x\in Ch_{B\setminus\{b\}}(i)} (bx)\otimes(\widehat{S(b)x}).
		\end{equation}
		
		(b) Deg=0. One particular derivation (the apparent asymmetry of $\theta_1$ and $\theta_2$ is due to SuperLie's aesthetic criteria)
		\begin{equation}\label{seeFn}
			D_{\theta}^{(0)} = \sum_{x\in\tilde{C}(\theta_1)} \theta_1x\otimes(\widehat{\theta_1x}).
		\end{equation}
		
		(c) Deg=0. Another particular derivation --- the Euler operator
		\begin{equation}\label{Eul}
		D_{0}^{(0)} = \sum_{x\in O} x \otimes(\widehat{x}).
		\end{equation}
	%\end{enumerate}

	(d) Deg = $2n-2$. See eq.~\eqref{2n-2}.

%\end{enumerate}

\subsection{$\mathfrak{out}(\fh_{I}^{(1)}(2n|0))$}\label{2a_0_I} The desuperization of the case $\mathfrak{out}(\fh_I^{(1)}(0|2n))$; the same cocycles.

\subsection{$\mathfrak{out}(\fh_{II}^{(1)}(2a|2b))$}\label{2a_2b_I}
Cocycles:

(a) Deg = $2a+2b-2$. See eq.~\eqref{2n-2}.
%\end{enumerate}

\subsection{$\mathfrak{out}(\mathfrak{h}^{(1)}(0|2n+1))$}\label{0_2b+1}
Notation for convenience:
\begin{enumerate}
	\item Basis: $B = \{\eta_1,\dots ,\eta_n; \theta; \xi_1,\dots ,\xi_n\}$ or $\tilde{B} = \{\eta_1,\dots ,\eta_n; \xi_1,\dots ,\xi_n\}$.
	\item Let $C(x)$ denote the set of all monomials each of which is a multiple of $x$;\\ let $\tilde{C}(x)$ denote the set of all monomials each not a multiple of $x$.
	\item The switch symbol: $S(\xi_i) = \eta_i$ and $S(\eta_i) = \xi_i$, whereas $S(\theta) = \theta$.
	\item Let $O$ denote the set of all monomials of odd degree.
\end{enumerate}
Cocycles:

(a) Deg=0. There are $2n$ equivalent derivations: of weight $2w$ for any $x\in \tilde{B}$ of weight $w$, 
		\begin{equation}\label{deg0x}
		D_x^{(0)}:=\sum_{y\in \tilde{C}(x)\bigcap \tilde{C}(S(x))} xy \otimes (\widehat{S(x)y}).
		\end{equation}

(b) Deg=0. One particular derivation --- the Euler operator, see eq.~\eqref{Eul}.
		%\begin{equation*}
		%D_{0}^{(0)} = \sum_{x\in O} x \otimes(\widehat{x}).
		%\end{equation*}

(c) Deg = $2n-1$: See eq.~\eqref{2n-2}.	
%\end{enumerate}
%\subsection{$\mathfrak{out}(\mathfrak{h}^{(1)}(0|2n+1))$}\label{2a+1_0}

\subsection{$\mathfrak{out}(\fh_{\Pi\Pi}^{(1)}(2a+1|2b+1))$}\label{2a+1_2b+1}
Cocycles:

(a) Deg=0. A particular derivation $D_{\theta}^{(0)}$, see eq.~\eqref{seeFn}
		%\begin{equation}
			%D_{\theta}^{(0)} = \sum_{x\in\tilde{C}(\theta_1)} \theta_1x\otimes(\widehat{\theta_1x}).
		%\end{equation}

(b) Deg=0. Another particular derivation --- the Euler operator, see eq.~\eqref{Eul}
		%\begin{equation*}
		%D_{0}^{(0)} = \sum_{x\in O} x \otimes(\widehat{x}).
		%\end{equation*}

(c) Deg = $2a+2b$. See eq.~\eqref{2n-2}.

\subsection{$\mathfrak{out}(\fh_{\Pi I}^{(1)}(2a|2b))$}\label{2a_2b_IPi}
Cocycles:

(a) Deg=0. There are $2b$ equivalent derivations, see eq.~\eqref{deg0x}.

(b) Deg=0. A particular derivation $D_{\theta}^{(0)}$, see eq.~\eqref{seeFn}
		%\begin{equation}
			%D_{\theta}^{(0)} = \sum_{x\in\tilde{C}(\theta_1)} \theta_1x\otimes(\widehat{\theta_1x}).
		%\end{equation}

(c) Deg=0. Another particular derivation --- the Euler operator, see eq.~\eqref{Eul}
		%\begin{equation*}
		%D_{0}^{(0)} = \sum_{x\in O} x \otimes(\widehat{x}).
		%\end{equation*}

(d) Deg = $2a+2b-2$. See eq.~\eqref{2n-2}.
%\end{enumerate}

\subsection{$\mathfrak{out}(\fh_{I\Pi}^{(1)}(2a|2b))$}\label{2a_2b_PiI} Same cocycles as for $\mathfrak{out}(\fh_{\Pi I}^{(1)}(2b|2a))$, but with $p,q$ and $\xi,\eta$ interchanged in all cocycles.

\subsection{$\mathfrak{out}(\fh_{\Pi\Pi}^{(1)}(2a|2b+1))$}\label{2a_2b+1PiPi} A partial desuperization of $\fh_{\Pi\Pi}^{(1)}(0|2n+1)$.

\subsection{$\mathfrak{out}(\fh_{\Pi\Pi}^{(1)}(2a+1|2b))$}\label{2a+1_2bPiPi} Same cocycles as for $\mathfrak{out}(\fh_{\Pi \Pi}^{(1)}(2b|2a+1))$, but with $p,q$ and $\xi,\eta$ interchanged in all cocycles.

\subsection{$\mathfrak{out}(\fle^{(1)}(n|n))$}\label{n_n_Le} Same cocycles as for $\mathfrak{out}(\fh_{\Pi}^{(1)}(2n|0))$ with $q,\pi$ replacing $q, p$; odd cocycles are marked in Subsection~\ref{0_2nPi}. Observe that $p(B)\equiv n\pmod 2$, so there are no $D_\od$-extensions for $n$ odd, and no $D_\ev$-extensions for $n$ even, since the center $c$ is always odd and $p(c)=p(B)+p(D)$.

\section{The general cases. Double extensions}\label{SgenDE}

Note: in all cases of $\fh$ series, $p(B_\fa)$ is congruent to the parity of the number of odd indeterminates; $p(B_\fa)\equiv a+\od$ for $\fle(a|a)$. Let $\bar X$ be the product of all indeterminates.

\subsection{The double extensions of $\fa=\fh^{(1)}(0|2n+1)$, $\fh_{\Pi\Pi}^{(1)}(2a|2b+1)$, $\fh_{\Pi\Pi}^{(1)}(2a+1|2b)$ for $a+b=n>2$}

\subsubsection{The $D_\ev$-extension}
The derivation $D_0^{(0)}$ does not preserve the NIS $B$. Indeed,
\[
B_{\fa}(D_0^{(0)}(\theta),\nfrac{\partial \bar X}{\partial \theta})=1 \text{ while } B_{\fa}(\theta,D_0^{(0)}(\nfrac{\partial \bar X}{\partial \theta}))=0.
\]
Now, the derivations $D_{x}^{(0)}$ preserve the bilinear form $B_{\fa}$, and the proof is similar to that for $\fh_\Pi^{(1)}(0|4)$ in Subsection~\ref{DEevB}; all respective double extensions are isomorphic. The double extension by means of the derivation $D_{x}^{(0)}$ for any $x\in\tilde B$, see Subsection~\ref{0_2b+1}, is a Lie superalgebra denoted by $\widetilde {\mathfrak{po}}(0|5)$ in \cite{BeB1} for $n=2$:
\begin{equation}\label{h05p}
\footnotesize
\renewcommand{\arraystretch}{1.4}
\begin{tabular}{|c|c|c|c|c|c|} \hline
Derivation & preserves $B_\fa$& $s_\fg(c)$ & Extension&Extension&Extension\\
\hline
$D_{x}^{(0)}$ & Yes & 0 & $\widetilde {\mathfrak{po}}(0|2n+1)$& $\widetilde {\mathfrak{po}}_{\Pi\Pi}(2a|2b+1)$& $\widetilde {\mathfrak{po}}_{\Pi\Pi}(2a+1|2b)$\\ \hline
$D_0^{(0)}$ & No & -- & -- & --& --\\ \hline
\end{tabular}
\end{equation}

\iffalse
\sssbegin{Lemma}\emph{([BeB1])} \emph{1)} Although $\sdim \fo\fo^{(1)}_{I\Pi}(4|4)=\sdim\widetilde{\mathfrak{po}}(0|5)=16|16$, see \emph{\cite{LeD}}, since $\dim\mathrm{H}^1(\widetilde{\mathfrak{po}}(0|5))=1$ and $\dim\mathrm{H}^1(\fo\fo^{(1)}_{I\Pi}(4|4))=0$, it follows that $\widetilde{\mathfrak{po}}(0|5)\not\simeq\fo\fo^{(1)}_{I\Pi}(4|4)$. \SJ{$\sdim\widetilde{\mathfrak{po}}(0|5) = 17|15$ because $\sdim\mathfrak{h}^{(1)}(0|5) = 15|15$ and we add an even center and an even derivation.}\DL{partly no: the derivation is odd}

 \emph{2)} \DL{Jin, we have to discuss the cases $n>5$ (unless you do it on your own}
\end{Lemma}\fi

%%%%%%%%%%%%%%%%%%%%%%%%%%%%%%%%%%%%%%%%%%%%%%%%%%%%%%%%%%%%%%
\subsubsection{The $D_\od$-extension}
%%%%%%%%%%%%%%%%%%%%%%%%%%%%%%%%%%%%%%%%%%%%%%%%%%%%%%%%%%%%%%
In a lexicographically ordered basis on $\fh^{(1)}(0|2n+1)$ the Gram matrix of NIS is $B_{\fa}=\antidiag(1,...,1)$.
The condition $B_{\fa}^TD^{(2n-1)}=D^{(2n-1)}B_{\fa}$ is easy to see. Besides, since $D^{(2n-1)}|_{\fg_{\ev}}=0$, it follows that
\[
B_{\fa}(D^{(2n-1)}(f),f)=0\text{~~ for any $f \in \fa_\ev$. }
\]
Now, since $0=(D^{(2n-1)})^2=\ad_{A}$, it follows that $A=0$ because $\fa$ has no center. We have, therefore, a parametric family of double extensions by means of $D^{(2n-1)}$, $q$, $A=0$ and $m\in \Kee$ (see Theorem \ref{MainTh3}). It is proved in \cite{BeB1}, that $\mathfrak{po}(0|5;m)$ for $m\neq 0$ is not isomorphic to $\mathfrak{po}(0|5;0):=\mathfrak{po}(0|5)$, whereas $\mathfrak{po}(0|5;m)\simeq \mathfrak{po}(0|5;\tilde m)$ for any pair $(m, \tilde m)$ such that $m\tilde m\neq 0$. The same arguments are true for any $n>2$. Table \eqref{h05} summarizes these results for any $m\in \Kee^\times$:
\begin{equation}\label{h05}
\footnotesize
\renewcommand{\arraystretch}{1.4}
\begin{tabular}{|c|c|c|c|c|c|} \hline
$D$ & $q(a)$ & $s_\fg(D)$ & Extension& Extension& Extension\\
\hline
$D^{(2n-1)}$ & yes&$0$ & $\mathfrak{po}(0|2n+1)$& ${\mathfrak{po}}_{\Pi\Pi}(2a|2b+1)$& $ {\mathfrak{po}}_{\Pi\Pi}(2a+1|2b)$\\ 
\hline
$D^{(2n-1)}$ & yes&$mc$& $\mathfrak{po}(0|2n+1;m)$& ${\mathfrak{po}}_{\Pi\Pi}(2a|2b+1;m)$& $ {\mathfrak{po}}_{\Pi\Pi}(2a+1|2b;m)$\\ \hline
\end{tabular}
\end{equation}

\subsection{The double extensions of $\fa=\fh_\Pi^{(1)}(0|2n)$ for $n>2$ and its desuperizations}
\begin{equation}\label{h04g}
\footnotesize
\renewcommand{\arraystretch}{1.4}
\begin{tabular}{|c|c|c|c|c|c|} \hline
$D$ & $q$&
$\fh_\Pi^{(1)}(0|2n)$&
$\fh_\Pi^{(1)}(2n|0)$&
$\fh_{\Pi\Pi}^{(1)}(2a|2b)$&$\fh_{\Pi\Pi}^{(1)}(2a+1|2b+1)$\\
\hline
$D_{b}^{(0)}$ & yes & $\widetilde{\mathfrak{po}}_\Pi(0|2n)$& $\widetilde{\mathfrak{po}}_\Pi(2n|0)$& $\widetilde{\mathfrak{po}}_{\Pi\Pi}(2a|2b)$&$\widetilde{\mathfrak{po}}_{\Pi\Pi}(2a+1|2b+1)$\\ 
\hline
$D_0^{(0)}$ & no& $-$& $-$& $-$& $-$\\ \hline
$D^{(2n-2)}$ & yes & $\mathfrak{po}_\Pi(0|2n)$& $\mathfrak{po}_\Pi(2n|0)$& $\mathfrak{po}_{\Pi\Pi}(2a|2b)$&$\mathfrak{po}_{\Pi\Pi}(2a+1|2b+1)$\\ 
\hline
\end{tabular}
\end{equation}

\subsection{The double extensions of $\fa=\fh_I^{(1)}(0|2n)$ for $n>2$ and its desuperizations} Clearly, the results in cases $I\Pi$ and $\Pi I$ are obtained from one another.
\begin{equation}\label{h04ii}
\footnotesize
\renewcommand{\arraystretch}{1.4}
\begin{tabular}{|c|c|c|c|c|c|} \hline
$D$ & $q$&
$\fh_I^{(1)}(0|2n)$&
$\fh_I^{(1)}(2n|0)$&
$\fh_{II}^{(1)}(2a|2b)$&
$\fh_{\Pi I}^{(1)}(2a|2b)$\\
\hline
$D_{0}^{(0)}$& yes & $\widehat{\mathfrak{po}}_I(0|2n)$& $\widehat{\mathfrak{po}}_{I}(2n|0)$
& $\widehat{\mathfrak{po}}_{II}(2a|2b)$& $\widehat{\mathfrak{po}}_{\Pi I}(2a|2b)$\\ \hline
$D_b^{(0)}$ & yes & $\widetilde{\mathfrak{po}}_I(0|2n)$& $\widetilde{\mathfrak{po}}_{I}(2n|0)$
& $\widetilde{\mathfrak{po}}_{II}(2a|2b)$& $\widetilde{\mathfrak{po}}_{\Pi I}(2a|2b)$\\ \hline
$D_\theta^{(0)}$ & $-$ & $-$& $-$& $-$& $-$\\ \hline
$D^{(2n-2)}$ & yes & $\mathfrak{po}_I(0|2n)$& $\mathfrak{po}_I(2n|0)$& $\mathfrak{po}_{II}(2a|2b)$& $\mathfrak{po}_{\Pi I}(2a|2b)$\\ \hline
\end{tabular}
\end{equation}

\subsection{The double extensions of $\fa=\fle^{(1)}(n|n)$}\label{le2n}

\begin{equation}\label{leGen}
\footnotesize
\renewcommand{\arraystretch}{1.4}
\begin{tabular}{|c|c|c||c|c|} \hline
$D$ & $A$&
$\fle^{(1)}(2n|2n)$& $A$&
$\fle^{(1)}(2n+1|2n+1)$\\
\hline
$D_{b}^{(0)}$ & $0$ & $\tilde\fb(2n|2n)$& $-$& $-$\\ 
\hline
$D_0^{(0)}$ & $-$ & $-$& $-$ & $-$\\ \hline
$D^{(4n-2)}$ & $0$ & $\mathfrak{b}(2n|2n)$& $0$&$\mathfrak{b}(2n+1|2n+1)$\\ \hline
\end{tabular}
\end{equation}

\sssbegin{Lemma} \emph{1)} The double extensions $\mathfrak{po}_B(a|b)$, $\widetilde{\mathfrak{po}}_B(a|b)$, and $\widehat{\mathfrak{po}}_B(a|b)$ are non-isomorphic to each other.

\emph{2)}  $\tilde{\mathfrak{b}}(2n|2n) \not\simeq \mathfrak{b}(2n|2n)$.
\end{Lemma}

\begin{proof} 1) We have $\rk D^{(2n-2)}=2n$ in $\mathfrak{po}_I(0|2n)$, and there is no such element in both $\widetilde{\mathfrak{po}}_I(0|2n)$ and $\widehat{\mathfrak{po}}_I(0|2n)$. The rank of $D_0^{(0)}$ in $\widehat{\mathfrak{po}}_I(0|2n)$ is $2^{2n-1}\neq 2n$ for $n>1$, and there is no such element in $\widetilde{\mathfrak{po}}_I(0|2n)$. 

The same type of arguments applies to other forms $B$ and superdimensions.

2) This statement is true because $\rk \ad_{p_i}= \rk \ad_{\pi_i} = 2^{4n-1}-1$ in $\tilde{\mathfrak{b}}(2n|2n)$ and no element in $\mathfrak{b}(2n|2n)$ has odd rank.
\end{proof}

\section{Remarks on double extensions. Cases of loop (super)algebras}\label{s2rem}

\ssbegin{Lemma} [Central extension must be nontrivial for DE to be indecomposable]\label{nontrCE} Let $\fg$ be a finite-dimensional Lie (super)algebra, and $B$ a NIS on it. Let $D\in\fder\ \fg$ preserve $B$, and the 2-cocycle $\omega(a,b):=B(Da,b)$ be trivial, i.e., there exists $\alpha\in\fg^*$ such that $\omega(a,b) = \alpha([a,b])$ for all $a,b\in\fg$. In this case, $D$ is an inner derivation, so the double extension is decomposable.
\end{Lemma}

\begin{proof}  Since $B$ is non-degenerate, there exists an $x\in \fg$ such that $\alpha(y) = B(x,y)$ for all $y\in\fg$. Then,  
 \[
 B(Da,b) = \alpha([a,b]) = B(x,[a,b]) = B([x,a],b) \text{~~for all $a,b\in\fg$}.
 \]
 For this $x$, we have $Da=[x,a]$ for all $a\in\fg$. \end{proof}  
 
\subsection{Remark on affine Kac--Moody (super)algebras constructed from loop (super)algebras in the case where the target (super)algebra has no outer derivations and non-trivial central extensions} Let $\fa$ be a simple finite-dimensional Lie (super)algebra
with a NIS $b$. Let $\fa^{\ell(1)}:=\fa\otimes \Cee[x^{-1}, x]$, where $x=\exp(i\varphi)$ for the angle parameter $\varphi$ on the circle, be the $\fa$-valued Lie (super)algebra of loops expandable into Laurent polynomials. It is easy to see that for any $n\in\Zee$, the bilinear form
\be\label{bilPne21}
\begin{array}{l}
B_n(f,g):=\Res b(f,g) x^n, \text{~~where}\\
\text{$\Res f(x)=$coeff. of $\frac{1}{x}$ in the Laurent series expansion of $f(x)$},
\end{array}
\ee
is a NIS on $\fa^{\ell(1)}$. Let $D_n=x^n\frac{d}{dx}$. The non-trivial central extension given by the cocycle with values in $\Cee c$
\be\label{hmm}
\omega(f,g):=\Res b(f,dg) =B_0(f, D_0(g))
\ee
make $\fg:=(\Cee c \ltimes \fa^{\ell(1)}) \ltimes \Cee D_1$ a double extension of $\fa^{\ell(1)}$ called \textit{affine Kac--Moody} Lie (super)algebra;  for the non-super cases, see \cite{K}. If $\fa$ has Cartan matrix, then $\fg$ also has Cartan matrix; the Dynkin diagram of $\fg$ is the extended Dynkin diagram of $\fa$.

The above description made us wonder: 

(A) In the case where $\fder\ \fa=\fa$ for simplicity, the space of outer derivations of 
$\fa^{\ell(1)}$ is $\fvect(1)=\fder\ \Cee[x^{-1}, x]$, What is so special in $D_1$ to be selected among the $D_n$, where $n\in\Zee$, for the role of outer derivation of $\fa^{\ell(1)}$ in construction of the double extension? 

(B) Why instead of the cocycle \eqref{bilPne2} with $D=D_1$, everybody uses eq.~\eqref{hmm} with $D=D_0$? 

Let $f=u\otimes x^s$ and $g=v\otimes x^t$ be $\fa$-valued functions with $u, v\in\fa$. Then,  condition \eqref{(1)} for $B_0$, see \eqref{bilPne21}, turns into
\[
(s+t)\Res b(u, v)x^{n+s+t-1} =0
\]
true if either $n+s+t\neq0$, or $s+t=0$. Therefore, $D_n$ preserves NIS $B_{-n}$, see \eqref{bilPne21}, on $\fa^{\ell(1)}$, and the cocycle proportional (with a non-zero factor) to \eqref{bilPne2}, see \eqref{hmm}, can be obtained as
\be\label{hm}
\omega(f,g):=B_n(f, D_{-n}g)\text{~~for any $n\in\Zee$}.
\ee
Since only $D_1$ is of degree 0, and only $D_0$ preserves $B_0$, we see that only $D_1$ and $B_0$ lead to a Cartan matrix of $\fg$. This answers the above questions (A) and (B) in this subsection.

\section{Are there too many double extensions of the loop algebras?}

Let us recall Vinberg's description of algebra, considered as a set with a structure (\cite{V}):  ``Actually, it is possible to consider any operation on any set.  However, only few algebraic structures are of real interest such as examples natural in the sense that they arose from the studies of the real world or the internal progress in mathematics.'' Having accepted this approach, we wonder how to modify the questions (C)--(E) below to make them ``of real interest''.

(C) How many equivalence classes of double extensions of the loop algebras are there?

Let $\fa$ be a simple Lie algebra. The central extensions obtained by means of derivations $D_n$ which preserve $B_{-n}$ are all equivalent for all $n$ since there is only one equivalence class. However, we can not prove that there is only one class of double extensions. For any $P\in \Cee[x, x^{-1}]$, consider a derivation of the loop algebras by setting $D_P := P\del$. For any \textbf{series} $R\in \Cee[[x, x^{-1}]]$, consider the form $B_R(f,g) = \Res Rfg$. It seems that $D_P$ preserves $B_R$ if and only if $PR\in\Cee^\times$ and all central extensions corresponding to various polynomials $P$ are proportional to each other. 

We have, however, no idea which of the numerous double extensions constructed from pairs $(D_P, B_{P^{-1}})$ are isomorphic to each other  (the case of the simple Lie \textbf{super}algebra $\fa$ having several outer derivations is even more complicated). Perhaps, this question C), as well as its extension to the twisted loops considered in the Subsection below, are not reasonable? 

\subsection{The cases where simple target superalgebra $\fa$ has both non-trivial central extension and outer derivation} Here is the list of all such algebras $\fa$. We consider two questions for the target algebras $\fa$ admitting a double extension:

D) Recall that if $\psi$ is an order $m$ automorphism of $\fg$, and $\fg_{\bar k}$, where $0\leq \bar k\leq m-1$, are eigenspaces of $\psi$ with eigenvalue $\exp(\nfrac{2\pi k i}{m})$, where $i=\sqrt{-1}$, then 
\[
\fg_\psi^{\ell(m)}:=\oplus_{0\leq \bar k\leq m-1,\ j\in\Zee}\ \fg_{\bar k}t^{k+mj}. 
\]
What can we say about double extensions of $\fa^{\ell(m)}$? For the list of automorphism of simple Lie superalgebras $\fa$, see tables in Section~\ref{Stables}.

E) What can we say about double extensions of $\fa^{\ell(m)}$ for $m>1$ if $\fa$ is a Lie superalgebra admitting a double extension, in particular, has a NIS? 

Perhaps, there is no need to separate  questions D) and E) and we are just being misguided by the fact that for $\psi=\Pi$ and 
$\psi=\Pi\circ(-\st)$ of Lie superalgebra $\fa=\fpsl(n|n)$, whose double extension has a Cartan matrix,  the twisted loop algebra $\fa^{\ell(m)}$ does not have double extension with Cartan matrix, and the other way round: the double extension of $\fb=\fpsq(n)$ has no Cartan matrix, whereas a double extension of $\fb^{\ell(2)}$ for $\psi=\delta_{-1}$ has a Cartan matrix. 

See also the case of $\psi=\per$ of  $\fa=\fosp_{\eps_3}(4|2)$ which has a Cartan matrix whereas $\fa^{\ell(3)}$ does not  have double extension with Cartan matrix.

\underline{$\fpsl(2|2)$}, see \cite{BGL}. There are three central extensions (the corresponding cocycles form the Lie algebra $\fsl(2)$)
\[
0\tto \fsl(2)\tto \fosp(4|2;0)\tto \fpsl(2|2)\tto 0.
\]
There are three outer derivations (the corresponding cocycles form the Lie algebra $\fsl(2)$)
\[
0\tto \fpsl(2|2)\tto \fosp(4|2;-1)\tto \fsl(2)\tto 0.
\]
To these three outer derivations correspond three double extensions of $\fpsl(2|2)\simeq\fh^{(1)}(0|4)$, see the 3rd column in table \eqref{h04}, where the ground field is different, but this is irrelevant. Two of these double extensions ---  $\fgl(2|2)$ and $\fpo(0|4)$ --- are expected due to the isomorphism $\fpsl(2|2)\simeq\fh^{(1)}(0|4)$, the third one, denoted $\widetilde{\fpo}(0|4)$, was found in \cite{BeB1}.

\underline{$\fpsl(n|n)$ with $n>2$}. There is one outer derivation, one central extension,  and one double extension: $\fgl(n|n)$.

\underline{$\fpsq(n)$ with $n>2$}. There is one outer derivation, one central extension,  and one double extension: $\fq(n)$.

\underline{$\fh^{(1)}(0|n)$ with $n>4$}. There are two outer derivations, one central extension,  and one double extension: $\fpo(0|n)$.

Therefore, we have the double extensions of loops with values in the above-listed double extensions:

 --- $\fgl(2|2)^{\ell(1)}$,  $\fpo(0|4)^{\ell(1)}$, and  $\widetilde{\fpo}(0|4)^{\ell(1)}$ are DEs of $\fpsl(2|2)^{\ell(1)}$
 
 ---  $\fgl(n|n)^{\ell(1)}$ for $n>2$ is a DE of  $\fpsl(n|n)^{\ell(1)}$
 
--- $\fq(n)^{\ell(1)}$ for $n>2$ is a DE of  $\fpsq(n)^{\ell(1)}$,

--- $\fpo(0|n)^{\ell(1)}$ for $n>4$ ais a DE of   $(\fh(0|n)^{(1)})^{\ell(1)}$ for $n>4$.

But there are also \textbf{other} double extensions of $\fpsl(2|2)^{\ell(1)}$,   as well as $\fpsl(n|n)^{\ell(1)}$ and $\fpsq(n)^{\ell(1)}$ for $n>2$ and $(\fh(0|n)^{(1)})^{\ell(1)}$ for $n>4$, determined by $\frac{d}{dx}$ as in eq.~\eqref{hmm} or $x\frac{d}{dx}$ as in eq.~\eqref{hm}.

There are also double extensions of twisted loops $\fg_\psi^{\ell(m)}$, with NIS on $\fg_\psi^{\ell(m)}$ induced by NIS on $\fg^{\ell(1)}$, see \cite{BKLS}.

%\vfill
%\newpage

\section{Tables} \label{Stables}\footnotesize

In Tables~\eqref{*a} and~\eqref{**a}, let $\Ad_A(X):=A^{-1}XA$ for any $X\in\fgl(m|n)$ and some even invertible $A\in\fgl(m|n)$. The automorphisms $\Pi$ and $\Pi\circ(-\st)$ of $\fpsl(n|n)$ are the ones induced from their namesake automorphisms of $\fgl(n|n)$.

\begin{table}[ht]\centering
{%
\caption{Finite order automorphisms of simple finite-dimensional Lie superalgebras $\fg$ over $\Cee$ and of their double extensions (\cite{Se})}\label{td4.0}\nopagebreak\tiny
\be\label{*a} \tabcolsep=3.5pt
 \begin{tabular}{|l|l|}
 \hlx{hv}
 $\per\in\Aut \fosp_{\eps_3} (4|2), \text{ where }
\eps_3:= \fnfrac{1}{2}(-1\pm i\sqrt{3})$&
 $\per(a,u)=((a_3,a_1,a_2), \eps  u_3\otimes u_1\otimes
 u_2)$\\
 \hlx{vhv}
  $d_{23}\in\Aut \fosp_\lambda(4|2),\text{ where }\RE\lambda=-\frac12$&
 $d_{23}(a,u)=((a_1,a_3,a_2), u_1\otimes u_3\otimes u_2)$ for any \\[1pt]
&
 $\arraycolsep=0pt\begin{array}{l}
(a_1,a_2,a_3)\in\fsl(V_1)\oplus\fsl(V_2)\oplus\fsl(V_3)= \fosp_\alpha(4|2)_{\ev}\\
 \text{and}\ u_1\otimes u_2\otimes u_3\in V_1\otimes V_2\otimes V_3=\fosp_\alpha(4|2)_{\od}
\end{array}$
 \\
 \hlx{vhv}
 $A\in\Aut\fpo(0|2n)$, \text{~~ where~~}$\fpo(0|2n)\simeq\Cee[\theta_1, \dots, \theta_{2n}]$&
 $A(\theta_i)=(-1)^{\delta_{1i}}\theta_i$,\\

 $B\in\Aut\fpo(0|2n)$&
$B(\theta_i)=\theta_i + \partial_{\theta_i}(\theta_1\ldots\theta_{2n})$,\\  %$B_t=\exp(\theta_1\ldots\theta_{2n}t)$\\
 \hlx{vhv}
 $\delta_\lambda\in \begin{cases} \Aut\fvect(0|n), \text{~~ where~~
 }\fvect(0|n)=\fder\Cee[\theta_1, \dots, \theta_n]\\ \Aut\fgl(n|m),\quad \text{where
 }\lambda\in\Cee^\times\end{cases}$ &
 $\arraycolsep=0pt
 \begin{array}{l}
 \delta_\lambda(\theta_i)=\lambda\theta_i\text{~~ for all~~
 }i\\
 \delta_\lambda\begin{pmatrix}A&B\\C&D\end{pmatrix}=\begin{pmatrix}A&\lambda B\\\lambda^{-1}C&D\end{pmatrix}\end{array}$\\
 \hlx{vhv}
 $\Pty\in  \Aut\fg$ &
 $\Pty(x)=(-1)^{p(x)}$ for any $x\in\fg$\\
\hlx{vhv}
 $\Ad_{J_{k,2n}(A)}\in \Aut\fosp(k|2n) \text{ where }  J_{k,2n}(A)=\diag(A,1_{2n})$&
 $\arraycolsep=0pt
 \begin{array}{l}
\text{for any } A\in\text{O}(k) %\\
\text{ such that }   \det A=-1,\  AA^t=1
 \end{array}$\\
 \hlx{vhv}
 $-\st\in\Aut\fgl(n|m)$&
 $-\st\begin{pmatrix}A&B\\C&D\end{pmatrix}=\begin{pmatrix}-A^t&C^t\\ -B^t&-D^t\end{pmatrix}$\\
  \hlx{vhv}
  $\Pi\in\Aut\fgl(n|n)$& $\Pi\begin{pmatrix}A&B\\ C&D\end{pmatrix}=\begin{pmatrix}D&C\\ B&A\end{pmatrix}$\\
 \hlx{vhv}
 $q\in\Aut\fq(n)$& $q: (A, B)\longmapsto(-A^t, i B^t)$\\
 %$q=-\st\circ\delta_{\sqrt{-1}}$\\
 \hlx{vh}
\end{tabular}
\ee
}
\end{table}
\[
\begin{minipage}{\textwidth}
\hskip 1.5em 
Set $T_{m,n}:=\diag(-1, 1_{2m-1}, 1_{2n})$ for $m>0$; let $J_{2n}=\antidiag(1_n, -1_n)$, set $B_{m,2n}:=\diag(1_m, J_{2n})$.\\
\hskip 1.5em The components $\fg_{\bar i}$ are described as $\fg_{\bar 0}$-modules. The~$\fosp(m|2n)=\fosp_B(V)$-module $S^2(V)$ is redu\-cible: $S^2(V)=S_0^2(V)\oplus \Cee B$.\\  
\hskip 1.5em Let $\eps_3:= \fnfrac{1}{2}(-1\pm i\sqrt{3})$. Automorphisms of $\fpsl(n|n)$ are induced by indicated automorphisms of $\fgl(n|n)$.\\
\hskip 1.5em In Table \ref{**a}, line 1, $\fg$ is any simple finite-dimensional Lie superalgebra.\\
\hskip 1.5em In Table \ref{**a}, line $8$, let $L^k$ be the irreducible $\fosp(1|2)$-module with highest weight $k$ and even highest weight vector. \\
\hskip 1.5em For $\fg_{\bar 0}\subset\fgl(V)$, we denote by $\id$ the tautological $\fg_\ev$-module $V$; let $\ad$ be the adjoint $\fg_{\bar 0}$-module. 
\end{minipage}
\]

{%\MathSkip{.5}%
\begin{table}[ht]\centering
{\footnotesize%
\caption{Lie superalgebras $\fg^{\ell(m)}_\psi$ and components
$\fg_i$ for $0\leq \bar i\leq m-1$ (after \cite{Se,LSS})} \label{t-add-3.1}
\vskip2mm
\be\label{**a}
 \tabcolsep=3.1pt
\begin{tabular}{|c|c|c|c|c|c|c|}
 \hlx{hv}
& $\fg$ & $\psi$ & $\fg_{\bar 0}$ & $\fg_{\bar 1}$ & $\fg_{\bar 2}$ & $\fg_{\bar 3}$\\
 \hlx{vhv}
$1$& $\fg$ & $\id$ & $\fg$ & $-$ & $-$ & $-$\\
 \hlx{vhhv}
$2$&  $\fsl(2m|2n)$ & $\Ad_{B_{2m,2n}}\circ(-\st)$ & $\fosp(2m|2n)$ &
 $S_0^2(\id)$ & $-$ & $-$\\
 \hlx{vhv}
$3$& $\fsl(2m+\mkern-2mu1|2n)$ & $\Ad_{B_{2m+\mkern-2mu1,2n}}\circ(-\st)$
 & $\fosp(2m+\mkern-2mu1|2n)$ & $S_0^2(\id)$&$\Pi(\id)$   & $\Pi(\id)$\\
% \hlx{vhv}
% $\fpsl(2n+\mkern-2mu1|2n+\mkern-2mu1)$ & $\Ad_{B_{2n+\mkern-2mu1,2n,1}}
% \circ(-\st)$ & $\fosp(2n+\mkern-2mu1|2n)$
% & $\Pi(\id)$ & $S^2(\id)/\Cee\langle\od\rangle$ & $\Pi(\id)$\\
 \hlx{vhv}
$4$&  $\fpsl(2n|2n)$ & $\Ad_{B_{2n,2n}}\circ(-\st)$ & $\fosp(2n|2n)$
 & $S_0^2(\id)$ & $-$ & $-$\\
 \hlx{vhv}
$5$&  $\fpsl(n|n)\text{  for  } n>2$ & $\Pi$ & $\fp\fq(n)$ & $\ad^*$ & $-$ & $-$\\
 \hlx{vhv}
$6$&  $\fpsl(n|n)\text{  for  } n> 2$ & $\Pi\circ(-\st)$ & $\fspe(n)$ &
 $\ad^*$ & $-$ & $-$\\
 \hlx{vhv}
$7$&  $\fosp(2m|2n)$ & $\Ad_{T_{m,n}}$ & $\fosp(2m-\mkern-2mu1|2n)$ & $\id$ & $-$ & $-$\\
 \hlx{vhv}
$8$&  $\fosp_{\eps_3}(4|2)$ & $\per$ & $\fosp(1|2)$ &$\ad=L^2$
 &$\Pi(L^3)$& $-$\\
 \hlx{vhv}
 $9$&  $\fh^{(1)}(2n)\text{  for  } n>2$ & induced by $A\in\Aut\fpo(0|2n)$ & $\fh(2n-\mkern-2mu1)$ & $\ad^*$ & $-$ & $-$\\
\hlx{vhhv}
 $10$&  $\fpsq(n) \text{  for  } n>2$ &induced by  $q\in\Aut\fq(n)$ & $\fo(n)$ & $\Pi(S_0^2(\id))$ &
 $S_0^2(\id)$ & $\Pi(E^2(\id))$\\
 \hlx{vhv}
$11$&  $\fpsq(n) \text{  for  } n>2$ & induced by  $\delta_{-\mkern-2mu1}\in\Aut\fq(n$ & $\fsl(n)$ & $\Pi(\ad)$ & $-$ & $-$\\
\hlx{vh}
\end{tabular}
\ee} \vspace{-3mm}
\end{table}
}

\normalsize

\vfill
\newpage

\textbf{Acknowledgements}. 
For the possibility to perform the difficult computations of this research we are grateful to M.~Al Barwani, Director of the High Performance Computing resources at New York University Abu Dhabi. We are thankful to J.~Bernstein, P.~Grozman, A.~Krutov, A.~Lebedev, and I.~Shchepochkina for helpful advice.
S.B. and D.L. were partly supported by the grant AD 065 NYUAD.

%%%%%

\end{document}